\documentclass[12pt,twoside,a4paper]{article}
\usepackage[english]{babel}
\usepackage{amssymb}
\usepackage[T1]{fontenc}
\usepackage{amsmath,epsf}
\newtheorem{Th}{Theorem}
\newtheorem{lemma}{Lemma}

\newcommand{\K}{\text{K}}
\newcommand{\E}{\mathbb{E}}
\newcommand{\F}{\mathcal{F}}
\newcommand{\R}{\mathbb{R}}
\newcommand{\Z}{\mathbb{Z}}
\newcommand{\N}{\mathbb{N}}

\newcommand{\B}{\mathcal{B}}

\newcommand{\U}{\mathcal{U}}
\newcommand{\V}{\mathcal{V}}

\renewcommand{\P}{\mathbb{P}}

\newcounter{tictac}
\newenvironment{fleuveA}{
   \begin{list}{($\textbf{A\arabic{tictac}}$) }{\usecounter{tictac}
\leftmargin 1cm\labelwidth 2em}}{\end{list}}

\def\1{\,\rlap{\mbox{\small\rm 1}}\kern.15em 1}

\def\build#1_#2^#3{\mathrel{\mathop{\kern 0pt#1}\limits_{#2}^{#3}}}
\def\tend#1#2{\build\hbox to 12mm{\rightarrowfill}_{#1\rightarrow #2}^{a.s.}}

\def\converge#1#2#3{\build\hbox to
15mm{\rightarrowfill}_{#1\rightarrow #2}^{\hbox{\scriptsize #3}}}
\begin{document}
\title{Asymptotic normality of the Parzen-Rosenblatt density estimator for strongly mixing random fields}
\author{Mohamed EL MACHKOURI,\\
Laboratoire de Math\'ematiques Rapha\"el Salem\\
UMR CNRS 6085, Universit\'e de Rouen\\
mohamed.elmachkouri@univ-rouen.fr}
\maketitle

{\renewcommand\abstractname{Abstract}
\begin{abstract}
\baselineskip=18pt 
We prove the asymptotic normality of the kernel density estimator (introduced by Rosenblatt (1956) and Parzen (1962)) in the context of stationary 
strongly mixing random fields. Our approach is based on the Lindeberg's method rather than on Bernstein's small-block-large-block technique and 
coupling arguments widely used in previous works on nonparametric estimation for spatial processes (see \cite{Carbon-Hallin-Tran}, 
\cite{Carbon-Tran-Wu}, \cite{Hallin-Lu-Tran-2001}, \cite{Tran}). Our method allows us to consider only minimal conditions 
on the bandwidth parameter and provides a simple criterion on the (non-uniform) strong mixing coefficients which do not depend on the bandwith. 
{\em Short title:} Kernel density estimation for random fields.
\end{abstract}
\thispagestyle{empty}
\baselineskip=18pt
\section{Introduction and main result}
The kernel density estimator introduced by Rosenblatt \cite{Ros} and Parzen \cite{Parzen1962} has received considerable 
attention in nonparametric estimation of probability densities for time series. If $(X_i)_{i\in\Z}$ is a stationary 
sequence of real random variables with a marginal density $f$ then the kernel density estimator of $f$ is defined for any positive integer $n$ and 
any $x$ in $\R$ by
$$
f_n(x)=\frac{1}{nb_n}\sum_{i=1}^n\K \left(\frac{x-X_i}{b_n}\right)
$$
where $\K$ is a probability kernel and the bandwidth $b_n$ is a parameter which converges slowly to zero such that $nb_n$ goes to infinity 
(the bandwidth determines the amount of smoothness of the estimator). For small $b_n$ we get a very rough estimate and for large $b_n$ a smooth estimate. 
The literature dealing with the asymptotic properties of $f_n$ when the observations $(X_i)_{i\in\Z}$ are independent 
is very extensive (see Silverman \cite{Silverman1986}). In particular, Parzen \cite{Parzen1962} proved that when $(X_i)_{i\in\Z}$ 
is i.i.d. and the bandwidth $b_n$ goes to zero such that $nb_n$ goes to infinity then $(nb_n)^{1/2}(f_n(x_0)-\E f_n(x_0))$ 
converges in distribution to the normal law with zero mean and variance $f(x_0)\int_{\R} \K ^2(t)dt$. This result was recently extended by Wu \cite{Wu-Mielniczuk} for 
causal linear processes with i.i.d. innovations under the same conditions on the bandwidth. Bosq, Merlev\`ede and Peligrad \cite{Bosq-Merlevede-Peligrad1999} established 
a central limit theorem for the kernel density estimator $f_n$ when the sequence $(X_i)_{i\in\Z}$ is assumed to be strongly mixing but the bandwith parameter $b_n$ 
is assumed to satisfy $b_n\geq C n^{-1/3}\log n$ (for some positive constant $C$) which is stronger than the bandwith parameter assumption in \cite{Parzen1962} and \cite{Wu-Mielniczuk}. 
In this paper, we are going to establish Parzen's central limit theorem (see Theorem $\ref{convergence-loi}$) for random variables which show spatial interaction 
(random fields). The problem is not trivial since $\Z^d$ do not have a natural ordering for $d\geq 2$ and consequently classical techniques available for one-dimensional 
processes do not extend to random fields. In particular, the martingale-difference method (Wu \cite{Wu-Mielniczuk}) for time series seems to be difficult 
to apply in the spatial context. Over the last few years nonparametric estimation for random fields (or spatial processes) was given increasing attention stimulated by a 
growing demand from applied research areas (see Guyon \cite{Guyon}). In fact, spatial data arise in various areas of 
research including econometrics, image analysis, meteorology, geostatistics... Key references on nonparametric 
estimation for random fields are Biau \cite{Biau-2003}, Carbon et al. \cite{Carbon-Hallin-Tran}, 
Carbon et al. \cite{Carbon-Tran-Wu}, Hallin et al. \cite{Hallin-Lu-Tran-2001}, \cite{Hallin-Lu-Tran-2004a}, 
Tran \cite{Tran}, Tran and Yakowitz \cite{Tran-Yakowitz} and Yao \cite{Yao} who have investigated nonparametric 
density estimation for random fields and Biau and Cadre \cite{Biau-Cadre}, 
El Machkouri \cite{Elmachkouri-SISP}, El Machkouri and Stoica \cite{Elmachkouri-Stoica}, Hallin et al. \cite{Hallin-Lu-Tran-2004b} 
and Lu and Chen \cite{Lu-Chen-2002}, \cite{Lu-Chen-2004} who have studied spatial prediction and spatial regression estimation.\\
Let $d$ be a positive integer and let $(X_i)_{i\in\Z^d}$ be a field of identically distributed real random variables with a 
marginal density $f$. Given two $\sigma$-algebras $\U$ and $\V$ of $\F$, different measures
of their dependence have been considered in the literature.
We are interested by one of them. The $\alpha$-mixing coefficient has
been introduced by Rosenblatt \cite{Ros} defined by
$$
\alpha(\U,\V)=\sup\{\vert\P(A\cap B)-\P(A)\P(B)\vert\, ,\,A\in\U,\,B\in\V\}.
$$
In the sequel, we consider the (non-uniform) strong mixing coefficients $\alpha_{1,\infty}(n)$ defined for each positive integer $n$ by
$$
\alpha_{1,\infty}(n)=\sup\,\{\alpha(\sigma(X_k),\F_{\Gamma}),\,
k\in\Z^d,\, \Gamma\subset\Z^d,\,\rho(\Gamma,\{k\})\geq n\},
$$
where $\F_{\Gamma}=\sigma(X_i\,;\,i\in\Gamma)$ and the distance $\rho$ is defined for any subsets $\Gamma_1$ and $\Gamma_2$ of $\Z^d$ by 
$\rho(\Gamma_{1},\Gamma_{2})=\min\{\vert
i-j\vert,\,i\in\Gamma_{1},\,j\in\Gamma_{2}\}$ with $\vert i-j\vert=\max_{1\leq s\leq d}\vert i_s-j_s\vert$ for any $i$ and $j$ in $\Z^d$. 
We say that the random field $(X_i)_{i\in\Z^d}$ is strongly mixing if 
$\lim_{n\to+\infty}\alpha_{1,\infty}(n)=0$. The class of mixing random fields in the above sense is very large (one can refer to Guyon \cite{Guyon} or 
Doukhan \cite{Doukhan} for examples) and we recall that Dedecker \cite{JD-tcl} obtained a central limit theorem for the stationary 
random field $(X_i)_{i\in\Z^d}$ provided that $X_0$ has zero mean and finite variance and
$$
\sum_{k\in\Z^d}\int_{0}^{\alpha_{1,\infty}(\vert k\vert)}Q^2_{X_0}(u)du<+\infty
$$
where $Q_{X_0}$ is the quantile 
function defined for any $u$ in $[0,1]$ by 
$$
Q_{X_0}(u)=\inf\{t\geq 0\,;\,\P(\vert X_0\vert>t)\leq u\}.
$$ 
We consider the density estimator of $f$ defined for any positive integer $n$ and any $x$ in $\R$ by
$$
f_n(x)=\frac{1}{n^db_n}\sum_{i\in\Lambda_n}\K\left(\frac{x-X_i}{b_n}\right)
$$
where $b_n$ is the bandwidth parameter, $\Lambda_n$ denotes the set $\{1,...,n\}^d$ and $\K $ is a probability kernel. Our aim is to provide a sufficient 
condition on the strong mixing coefficients $\alpha_{1,\infty}(n)$ for 
$(n^db_n)^{1/2}(f_n(x_i)-\E  f_n(x_i))_{1\leq i\leq k},\,(x_i)_{1\leq i\leq k}\in\R^k,\,k\in\N^{\ast},$ to 
converge in law to a multivariate normal distribution (Theorem $\ref{convergence-loi}$) under minimal conditions on the bandwidths 
(that is $b_n$ goes to zero and $n^db_n$ goes to infinity).\\
\\
We consider the following assumptions:
\begin{fleuveA}
\item The marginal probability distribution of each $X_k$ is absolutely continuous with continuous positive density function $f$.
\item The joint probability distribution of each $(X_0,X_k)$ is absolutely continuous with continuous joint density $f_{0,k}$.
\item $\K$ is a probability kernel with compact support and $\int_{\R} \K ^2(u)\,du<\infty$.
\item The bandwidth $b_n$ converges to zero and $n^db_n$ goes to infinity.
\end{fleuveA}
Our main result is the following.
\begin{Th}\label{convergence-loi}
Assume that $\textbf{\emph{(A1)}}$, $\textbf{\emph{(A2)}}$, $\textbf{\emph{(A3)}}$ and $\textbf{\emph{(A4)}}$ hold and 
\begin{equation}\label{mixing-condition}
\sum_{m=1}^{+\infty}m^{2d-1}\,\alpha_{1,\infty}(m)<+\infty.
\end{equation}
Then for any positive integer $k$ and any distinct points $x_1,...,x_k$ in $\R$,
\begin{equation}\label{limit}
(n^dh_{n})^{1/2}
\left(\begin{array}{c}
       f_{n}(x_1)-\E f_n(x_1)\\
       \vdots\\
       f_{n}(x_k)-\E f_n(x_k)
       \end{array} \right)
\converge{n}{+\infty}{\textrm{$\mathcal{L}$}}
\mathcal{N}\left(0,V\right)
\end{equation}
where $V$ is a diagonal matrix with diagonal elements $v_{ii}=f(x_i)\int_{\R}\K^2(u)du$.
\end{Th}
\textbf{Remark 1.} A replacement of $\E f_n(x_i)$ by $f(x_i)$ for any $1\leq i\leq k$ in ($\ref{limit}$) 
is a classical problem in density estimation theory. For example, if $f$ is assumed to be Lipschitz and 
if $\int_{\R}\vert u\vert \vert \K (u)\vert du<\infty$ then $\vert \E f_n(x_i)-f(x_i)\vert=O(b_n)$ and thus the 
centering $\E f_n(x_i)$ may be changed to $f(x_i)$ without affecting the above result provided that $n^db_n^3$ converges to zero.\\
\\
\textbf{Remark 2.} Theorem $\ref{convergence-loi}$ is an extension of Theorem 3.1 by Bosq, Merlev\`ede and Peligrad \cite{Bosq-Merlevede-Peligrad1999}. 
In fact, using a different approach, the authors obtained the same result for $d=1$ with an additional condition on the bandwith 
parameter: there exists a positive constant $C$ such that $b_n\geq C\,n^{-1/3}\log n$.
\section{Proofs}
\textbf{Proof of Theorem $\textbf{\ref{convergence-loi}}$}. Without loss of generality, we consider only the case $k=2$ and we 
refer to $x_1$ and $x_2$ as $x$ and $y$ ($x\neq y$). Let $\lambda_1$ and $\lambda_2$ be two constants such that $\lambda_1^2+\lambda_2^2=1$ and 
denote 
$$
S_n=\lambda_1(n^db_n)^{1/2}(f_n(x)-\E  f_n(x))+\lambda_2(n^db_n)^{1/2}(f_n(y)-\E  f_n(y))
=\sum_{i\in\Lambda_n}\frac{\Delta_i}{n^{d/2}}
$$ 
where $\Delta_i=\lambda_1Z_{i}(x)+\lambda_2Z_{i}(y)$ and for any $z$ in $\R$,
$$
Z_{i}(z)=\frac{1}{\sqrt{b_n}}\left(\K \left(\frac{z-X_{i}}{b_n}\right)-\E \K \left(\frac{z-X_{i}}{b_n}\right)\right).
$$
We consider the notations
\begin{equation}\label{eta_sigma}
\eta=(\lambda_1^2f(x)+\lambda_2^2f(y))\sigma^2\quad\textrm{and}\quad\sigma^2=\int_{\R} \K ^2(u)du.
\end{equation}
The proof of the following technical result is postponed to the annex.
\begin{lemma}\label{lemme-technique} $\E(\Delta_0^2)$ converges to $\eta$ and $\E\vert\Delta_0\Delta_i\vert=O(b_n)$ for any $i$ in $\Z^d\backslash\{0\}$.
\end{lemma}
In order to prove the convergence in distribution of $S_n$ to
$\sqrt{\eta}\tau_0$ where $\tau_0\sim\mathcal{N}(0,1)$, we are going to follow the Lindeberg's method used in the proof of 
the central limit theorem for stationary random fields by Dedecker \cite{JD-tcl}. Let us note that several previous asymptotic 
results for kernel density estimates in the context of spatial processes were established using the so-called Bernstein's 
small-block-large-block technique and coupling arguments which lead to restrictive conditions on the bandwith parameter (see for example \cite{Carbon-Hallin-Tran}, 
\cite{Carbon-Tran-Wu}, \cite{Hallin-Lu-Tran-2001}, \cite{Tran}). 
Our approach seems to be better since we obtain a central limit theorem when the bandwith satisfies only Assumption ($\textbf{A2}$).\\
\\
Let $\mu$ be the law of the stationary real random field
$(X_{k})_{k\in\Z^{d}}$ and consider the projection $\pi_0$
from $\R^{\Z^{d}}$ to $\R$ defined by $\pi_0(\omega)=\omega_{0}$ and
the family of translation operators $(T^{k})_{k\in\Z^{d}}$ from
$\R^{\Z^{d}}$ to $\R^{\Z^{d}}$ defined by
$(T^{k}(\omega))_{i}=\omega_{i+k}$ for any $k\in\Z^{d}$ and any
$\omega$ in $\R^{\Z^{d}}$. Denote by $\B$ the Borel
$\sigma$-algebra of $\R$. The random field $(\pi_0\circ
T^{k})_{k\in\Z^{d}}$ defined on the probability space
$(\R^{\Z^{d}}, \B^{\Z^{d}}, \mu)$ is stationary with the same law
as $(X_{k})_{k\in\Z^{d}}$, hence, without loss of
generality, one can suppose that $(\Omega, \F, \P)=(\R^{\Z^{d}},
\B^{\Z^{d}}, \mu)$ and $X_{k}=\pi_0\circ T^{k}$. On the lattice $\Z^{d}$ we define the
lexicographic order as follows: if $i=(i_{1},...,i_{d})$ and
$j=(j_{1},...,j_{d})$ are distinct elements of $\Z^{d}$, the
notation $i<_{\text{\text{lex}}}j$ means that either $i_{1}<j_{1}$ or for some
$p$ in $\{2,3,...,d\}$, $i_{p}<j_{p}$ and $i_{q}=j_{q}$ for $1\leq
q<p$. Let the sets
$\{V_{i}^{M}\,;\,i\in\Z^{d}\,,\,M\in\N^{\ast}\}$ be defined as
follows:
$$
V_{i}^{1}=\{j\in\Z^{d}\,;\,j<_{\text{lex}}i\},
$$
and for $M\geq 2$
$$
V_{i}^{M}=V_{i}^{1}\cap\{j\in\Z^{d}\,;\,\vert i-j\vert\geq
M\}\quad\textrm{where}\quad \vert i-j\vert=\max_{1\leq l\leq
d}\vert i_{l}-j_{l}\vert.
$$
For any subset $\Gamma$ of $\Z^{d}$ define
$\F_{\Gamma}=\sigma(X_{i}\,;\,i\in\Gamma)$ and set
$$
\E_{M}(X_{i})=\E(X_{i}\vert\F_{V_{i}^{M}}),\quad M\in \N^{\ast}.
$$
Let $g$ be a one to one map from $[1,M]\cap\N^{\ast}$ to a finite
subset of $\Z^d$ and $(\xi_i)_{i\in\Z^d}$ a real random field. For
all integers $k$ in $[1,M]$, we denote
$$
S_{g(k)}(\xi)=\sum_{i=1}^k \xi_{g(i)}\quad\textrm{and}\quad
S_{g(k)}^{c}(\xi)=\sum_{i=k}^M \xi_{g(i)}
$$
with the convention $S_{g(0)}(\xi)=S_{g(M+1)}^{c}(\xi)=0$. To
describe the set $\Lambda_{n}=\{1,...,n\}^d$, we define the one to
one map $g$ from $[1,n^d]\cap\N^{\ast}$ to $\Lambda_{n}$ by:
$g$ is the unique function such that $g(k)<_{\text{\text{lex}}}g(l)$ for $1\leq k<l\leq
n^d$. From now on, we consider a field $(\tau_{i})_{i\in\Z^d}$ of i.i.d. random variables independent of
$(X_{i})_{i\in\Z^d}$ such that $\tau_{0}$ has the standard normal
law $\mathcal{N}(0,1)$. We introduce the fields $Y$ and $\gamma$ defined for any $i$ in $\Z^d$ by 
$$
Y_{i}=\frac{\Delta_i}{n^{d/2}}\quad\textrm{and}\quad\gamma_{i}=\frac{\tau_i\sqrt{\eta}}{n^{d/2}}
$$ 
where $\eta$ is defined by ($\ref{eta_sigma}$).\\
\\
Let $h$ be any function from $\R$ to $\R$. For $0\leq k\leq l\leq n^d +1$, we introduce
$h_{k,l}(Y)=h(S_{g(k)}(Y)+S_{g(l)}^{c}(\gamma))$. With the above
convention we have that $h_{k,n^d+1}(Y)=h(S_{g(k)}(Y))$ and also
$h_{0,l}(Y)=h(S_{g(l)}^{c}(\gamma))$. In the sequel, we will often
write $h_{k,l}$ instead of $h_{k,l}(Y)$. We denote by
$B_{1}^4(\R)$ the unit ball of $C_{b}^4(\R)$: $h$ belongs to
$B_{1}^4(\R)$ if and only if it belongs to $C^4(\R)$ and satisfies
$\max_{0\leq i\leq 4}\|h^{(i)}\|_{\infty}\leq 1$.\\
\\
It suffices to prove that for all $h$ in $B_{1}^4(\R)$,
$$
\E\left(h\left(S_{g(n^d)}(Y)\right)\right)\converge{n}{+\infty}{}\E \left(h\left(\tau_0\sqrt{\eta}\right)\right).
$$
We use Lindeberg's decomposition:
$$
\E \left(h\left(S_{g(n^d)}(Y)\right)-h\left(\tau_{0}\sqrt{\eta}\right)\right)
=\sum_{k=1}^{n^d}\E \left(h_{k,k+1}-h_{k-1,k}\right).
$$
Now,
$$
h_{k,k+1}-h_{k-1,k}=h_{k,k+1}-h_{k-1,k+1}+h_{k-1,k+1}-h_{k-1,k}.
$$
Applying Taylor's formula we get that:
$$
h_{k,k+1}-h_{k-1,k+1}=Y_{g(k)}h_{k-1,k+1}^{'}+\frac{1}{2}Y_{g(k)}^{2}h_{k-1,k+1}^{''}+R_{k}
$$
and
$$
h_{k-1,k+1}-h_{k-1,k}=-\gamma_{g(k)}h_{k-1,k+1}^{'}-\frac{1}{2}\gamma_{g(k)}^{2}h_{k-1,k+1}^{''}+r_{k}
$$
where $\vert R_{k}\vert\leq Y_{g(k)}^2(1\wedge\vert Y_{g(k)}\vert)$ and $\vert r_{k}\vert\leq\gamma_{g(k)}^2(1\wedge\vert\gamma_{g(k)}\vert)$.
Since $(Y,\tau_{i})_{i\neq g(k)}$ is independent of
$\tau_{g(k)}$, it follows that
$$
\E \left(\gamma_{g(k)}h_{k-1,k+1}^{'}\right)=0\quad\textrm{and}\quad
\E \left(\gamma_{g(k)}^2h_{k-1,k+1}^{''}\right)=\E \left(\frac{\eta}{n^d}h_{k-1,k+1}^{''}\right)
$$
Hence, we obtain
\begin{align*}
\E \left(h(S_{g(n^d)}(Y))-h\left(\tau_0\sqrt{\eta}\right)\right)&=
\sum_{k=1}^{n^d}\E (Y_{g(k)}h_{k-1,k+1}^{'})\\
&\quad+\sum_{k=1}^{n^d}\E \left(\left(Y_{g(k)}^2-\frac{\eta}{n^d}\right)\frac{h_{k-1,k+1}^{''}}{2}\right)\\
&\quad+\sum_{k=1}^{n^d}\E \left(R_{k}+r_{k}\right).
\end{align*}
Let $1\leq k\leq n^d$ be fixed. Noting that $\Delta_0$ is bounded by $4\|\K\|_{\infty}/\sqrt{b_n}$ and applying Lemma $\ref{lemme-technique}$, 
we derive
$$
\E\vert R_k\vert\leq\frac{\E\vert\Delta_0\vert^3}{n^{3d/2}}=O\left(\frac{1}{(n^{3d}\,b_n)^{1/2}}\right)
$$
and
$$
\E\vert r_k\vert\leq\frac{\E\vert\gamma_0\vert^3}{n^{3d/2}}\leq\frac{\eta^{3/2}\E\vert\tau_0\vert^3}{n^{3d/2}}=O\left(\frac{1}{n^{3d/2}}\right).
$$
Consequently, we obtain
$$
\sum_{k=1}^{n^d}\E \left(\vert R_{k}\vert+\vert r_{k}\vert\right)=O\left(\frac{1}{(n^db_n)^{1/2}}+\frac{1}{n^{d/2}}\right)=o(1).
$$
Now, it is sufficient to show
\begin{equation}\label{equation1}
\lim_{n\to+\infty}\sum_{k=1}^{n^d}\left(\E (Y_{g(k)}h_{k-1,k+1}^{'})+\E \left(\left(Y_{g(k)}^2-\frac{\eta}{n^d}\right)\frac{h_{k-1,k+1}^{''}}{2}\right)\right)=0.
\end{equation}
First, we focus on
$\sum_{k=1}^{n^d}\E \left(Y_{g(k)}h_{k-1,k+1}^{'}\right)$. For all
$M$ in $\N^{\ast}$ and all integer $k$ in $[1,n^d]$, we define
$$
E_{k}^M=g([1,k]\cap\N^{\ast})\cap V_{g(k)}^M\quad\textrm{and}\quad
S_{g(k)}^M(Y)=\sum_{i\in E_{k}^M}Y_{i}.
$$
For any function $\Psi$ from $\R$ to $\R$, we define $\Psi_{k-1,l}^M=\Psi(S_{g(k)}^M(Y)+S_{g(l)}^c(\gamma))$ (we are going to 
apply this notation to the successive derivatives of the function $h$).\\
For any integer $n$, we define 
$$
m_n=\max\left\{\bigg[b_n^{\frac{-1}{2d}}\bigg],\left[\left(\frac{1}{b_n^2}\sum_{\vert i\vert>\left[b_n^{\frac{-1}{2d}}\right]}
\vert i\vert^d\,\alpha_{1,\infty}(\vert i\vert)\right)^{\frac{1}{2d}}\right]+1\right\}
$$
where $[\,.\,]$ denotes the integer part function. The following technical result will be proved in the annex. 
\begin{lemma}\label{lemme-conditions-on-m_n}
Under Assumption $\textbf{\emph{(A4)}}$ and the mixing condition $(\ref{mixing-condition})$, we have
\begin{equation}\label{conditions-on-m_n}
m_n^d\rightarrow\infty,\qquad m_n^db_n\rightarrow 0\qquad\textrm{and}\qquad
\frac{1}{m_n^db_n}\sum_{\vert i\vert>m_n}\vert i\vert^d\,\alpha_{1,\infty}(\vert i\vert)\rightarrow 0.
\end{equation}
\end{lemma}
Our aim is to show that
$$
\lim_{n\to+\infty}\sum_{k=1}^{n^d}\E \left(\left(Y_{g(k)}h_{k-1,k+1}^{'}-Y_{g(k)}\left(S_{g(k-1)}(Y)-S_{g(k)}^{m_n}(Y)\right)h_{k-1,k+1}^{''}\right)\right)=0.
$$
First, we use the decomposition
$$
Y_{g(k)}h_{k-1,k+1}^{'}=Y_{g(k)}h_{k-1,k+1}^{'m_n}+Y_{g(k)}\left(h_{k-1,k+1}^{'}-h_{k-1,k+1}^{'m_n}\right).
$$
We consider a one to one map $m$ from $[1,\vert E_{k}^{m_n}\vert]\cap\N^{\ast}$ to $E_{k}^{m_n}$ and such that $\vert
m(i)-g(k)\vert\leq\vert m(i-1)-g(k)\vert$. This choice of $m$
ensures that $S_{m(i)}(Y)$ and $S_{m(i-1)}(Y)$ are
$\F_{V_{g(k)}^{\vert m(i)-g(k)\vert}}$-measurable. The fact that
$\gamma$ is independent of $Y$ imply that
$$
\E \left(Y_{g(k)}h^{'}\left(S_{g(k+1)}^c(\gamma)\right)\right)=0.
$$
Therefore
\begin{equation}\label{equation_theta}
\left\vert \E \left(Y_{g(k)}h_{k-1,k+1}^{'{m_n}}\right)\right\vert=\left\vert\sum_{i=1}^{\vert
E_{k}^{m_n}\vert}\E \left(Y_{g(k)}\left(\theta_i-\theta_{i-1}\right)\right)\right\vert
\end{equation}
where $\theta_i=h^{'}\left(S_{m(i)}(Y)+S_{g(k+1)}^c(\gamma)\right)$.\\
\\
Since $S_{m(i)}(Y)$ and $S_{m(i-1)}(Y)$ are $\F_{V_{g(k)}^{\vert
m(i)-g(k)\vert}}$-measurable, we can take the conditional
expectation of $Y_{g(k)}$ with respect to $\F_{V_{g(k)}^{\vert
m(i)-g(k)\vert}}$ in the right hand side of ($\ref{equation_theta}$). On
the other hand the function $h^{'}$ is $1$-Lipschitz, hence
$$
\left\vert\theta_i-\theta_{i-1}\right\vert\leq\vert Y_{m(i)}\vert.
$$
Consequently,
$$
\left\vert
\E \left(Y_{g(k)}\left(\theta_i-\theta_{i-1}\right)\right)\right\vert\leq\E \vert Y_{m(i)}\E _{\vert m(i)-g(k)\vert}\left(Y_{g(k)}\right)\vert
$$
and
$$
\left\vert
\E \left(Y_{g(k)}h_{k-1,k+1}^{'m_n}\right)\right\vert\leq\sum_{i=1}^{\vert
E_{k}^{m_n}\vert}\E \vert Y_{m(i)}\E _{\vert m(i)-g(k)\vert}(Y_{g(k)})\vert.
$$
Hence,
\begin{align*}
\left\vert\sum_{k=1}^{n^d}\E \left(Y_{g(k)}h_{k-1,k+1}^{'m_n}\right)\right\vert
&\leq\frac{1}{n^d}\sum_{k=1}^{n^d}\sum_{i=1}^{\vert
E_{k}^{m_n}\vert} \E \vert\Delta_{m(i)}\E _{\vert
m(i)-g(k)\vert}(\Delta_{g(k)})\vert\\
&\leq \sum_{\vert j\vert \geq m_n}\|\Delta_j\E _{\vert j\vert}(\Delta_0)\|_1.
\end{align*}
For any $j$ in $\Z^d$, we have 
$$
\|\Delta_j\E _{\vert j\vert}(\Delta_0)\|_1
=\textrm{Cov}\left(\vert\Delta_j\vert\left(\mathbb{I}_{\E_{\vert j\vert}(\Delta_0)\geq 0}-\mathbb{I}_{\E_{\vert j\vert}(\Delta_0)<0}\right),\Delta_0\right).
$$
So, applying Rio's covariance inequality (cf. \cite{Rio1993}, Theorem 1.1), we obtain
$$
\|\Delta_j\E _{\vert j\vert}(\Delta_0)\|_1\leq 4\int_{0}^{\alpha_{1,\infty}(\vert j\vert)}Q_{\Delta_0}^2(u)du
$$
where $Q_{\Delta_0}$ is defined by $Q_{\Delta_0}(u)=\inf\{t\geq 0\,;\,\P(\vert \Delta_0\vert>t)\leq u\}$ for any $u$ in $[0,1]$. Since $\Delta_0$ is bounded by 
$4\|\K\|_{\infty}/\sqrt{b_n}$, we have
$$
Q_{\Delta_0}(u)\leq\frac{4\|\K\|_{\infty}}{\sqrt{b_n}}\qquad\textrm{and}\qquad
\|\Delta_j\E _{\vert j\vert}(\Delta_0)\|_1\leq\frac{64\|\K\|_{\infty}^2}{b_n}\,\alpha_{1,\infty}(\vert j\vert).
$$
Finally, we derive
\begin{align*}
\left\vert\sum_{k=1}^{n^d}\E \left(Y_{g(k)}h_{k-1,k+1}^{'m_n}\right)\right\vert
&\leq\frac{64\|K\|_{\infty}^2}{b_n}\sum_{\vert j\vert \geq m_n}\alpha_{1,\infty}(\vert j\vert)\\
&\leq\frac{64\|K\|_{\infty}^2}{m_n^db_n}\sum_{\vert j\vert \geq m_n}\vert j\vert^d\,\alpha_{1,\infty}(\vert j\vert)\\
&=o(1)\qquad\textrm{by ($\ref{conditions-on-m_n}$)}.
\end{align*}
Applying again Taylor's formula, it remains to consider
$$
Y_{g(k)}(h_{k-1,k+1}^{'}-h_{k-1,k+1}^{'m_n})=Y_{g(k)}(S_{g(k-1)}(Y)-S_{g(k)}^{m_n}(Y))h_{k-1,k+1}^{''}+R_{k}^{'},
$$
where $\vert R_{k}^{'}\vert\leq 2\vert Y_{g(k)}(S_{g(k-1)}(Y)-S_{g(k)}^{m_n}(Y))(1\wedge\vert S_{g(k-1)}(Y)-S_{g(k)}^{m_n}(Y)\vert)\vert$. 
Denoting $W_n=\{-m_n+1,...,m_n-1\}^d$ and $W_n^{\ast}=W_n\backslash\{0\}$, it follows that
\begin{align*}
\sum_{k=1}^{n^d}\E \vert R_{k}^{'}\vert 
&\leq 2\E \left(\vert\Delta_{0}\vert\left(\sum_{i\in W_n}\vert\Delta_{i}\vert\right)
\left(1\wedge\frac{1}{n^{d/2}}\sum_{i\in W_n}\vert\Delta_{i}\vert\right)\right)\\
&=2\E\left(\left(\Delta_0^2+\sum_{i\in W_n^{\ast}}\vert\Delta_0\Delta_i\vert\right)\left(1\wedge\frac{1}{n^{d/2}}\sum_{i\in W_n}\vert\Delta_i\vert\right)\right)\\
&\leq\frac{2}{n^{d/2}}\sum_{i\in W_n}\E(\Delta_0^2\vert\Delta_i\vert)+2\sum_{i\in W_n^{\ast}}\E\vert\Delta_0\Delta_i\vert\\
&\leq\frac{8\|\K\|_{\infty}}{(n^{d}b_n)^{1/2}}\sum_{i\in W_n}\E(\vert\Delta_0\Delta_i\vert)+2\sum_{i\in W_n^{\ast}}\E\vert\Delta_0\Delta_i\vert
\quad\textrm{since $\Delta_0\leq\frac{4\|K\|_\infty}{\sqrt{b_n}}$ a.s.}\\
&=\frac{8\E (\Delta_0^2)\|\K\|_{\infty}}{(n^{d}b_n)^{1/2}}+2\left(1+\frac{4\|\K\|_{\infty}}{(n^{d}b_n)^{1/2}}\right)\sum_{i\in W_n^{\ast}}\E(\vert\Delta_0\Delta_i\vert)\\
&=O\left(\frac{1}{(n^{d}b_n)^{1/2}}+m_n^db_n\left(1+\frac{1}{(n^{d}b_n)^{1/2}}\right)\right)\qquad\textrm{(by Lemma $\ref{lemme-technique}$)}\\
&=o(1)\quad\textrm{by ($\ref{conditions-on-m_n}$)}.
\end{align*}
So, we have shown that 
$$
\lim_{n\to +\infty}\sum_{k=1}^{n^d}\E \left(Y_{g(k)}h^{'}_{k-1,k+1}-Y_{g(k)}(S_{g(k-1)}-S_{g(k)}^{m_n})h^{''}_{k-1,k+1}\right)=0.
$$
In order to obtain ($\ref{equation1}$) it remains to control
$$
F_{0}=\E \left(\sum_{k=1}^{n^d}h_{k-1,k+1}^{''}\left(\frac{Y_{g(k)}^2}{2}+Y_{g(k)}\left(S_{g(k-1)}(Y)-S_{g(k)}^{m_n}(Y)\right)-
\frac{\eta}{2n^d}\right)\right).
$$
We consider the following sets:
$$
\Lambda_{n}^{m_n}=\{i\in\Lambda_{n}\,;\,\rho(\{i\},\partial\Lambda_{n})\geq m_n\}
\quad\textrm{and}\quad I_{n}^{m_n}=\{1\leq i\leq n^d\,;\,g(i)\in\Lambda_{n}^{m_n}\},
$$
and the function $\Psi$ from $\R^{\Z^d}$ to $\R$ such that
$$
\Psi(\Delta)=\Delta_{0}^2+\sum_{i\in V_{0}^1\cap W_n}2\Delta_{0}\Delta_{i}\quad\textrm{where $W_n=\{-m_n+1,...,m_n-1\}^d$}.
$$
For $1\leq k\leq n^d$, we set $D_{k}^{(n)}=\eta-\Psi\circ T^{g(k)}(\Delta)$. By definition of $\Psi$ and of the set
$I_{n}^{m_n}$, we have for any $k$ in $I_{n}^{m_n}$
$$
\Psi\circ T^{g(k)}(\Delta)=\Delta_{g(k)}^2+2\Delta_{g(k)}(S_{g(k-1)}(\Delta)-S_{g(k)}^{m_n}(\Delta)).
$$
Therefore for $k$ in $I_{n}^{m_n}$
$$
\frac{D_{k}^{(n)}}{n^d}=\frac{\eta}{n^d}-Y_{g(k)}^2-2Y_{g(k)}(S_{g(k-1)}(Y)-S_{g(k)}^{m_n}(Y)).
$$
Since $\lim_{n\to+\infty}n^{-d}\vert I_{n}^{m_n}\vert=1$, it remains to consider
$$
F_1=\left\vert\E\left(\frac{1}{n^d}\sum_{k=1}^{n^d}h_{k-1,k+1}^{''}D_{k}^{(n)}\right)\right\vert.
$$
Applying Lemma $\ref{lemme-technique}$, we have
\begin{align*}
F_1&\leq\left\vert\E\left(\frac{1}{n^d}\sum_{k=1}^{n^d}h_{k-1,k+1}^{''}(\Delta_{g(k)}^2-\E(\Delta_0^2))\right)\right\vert
+\vert\eta-\E(\Delta_0^2)\vert+2\sum_{j\in V_0^1\cap W_n}\E\vert\Delta_0\Delta_{j}\vert\\
&\leq\left\vert\E\left(\frac{1}{n^d}\sum_{k=1}^{n^d}h_{k-1,k+1}^{''}(\Delta_{g(k)}^2-\E(\Delta_0^2))\right)\right\vert
+o(1)+O(m_n^db_n),
\end{align*}
it suffices to prove that
$$
F_2=\left\vert\E\left(\frac{1}{n^d}\sum_{k=1}^{n^d}h_{k-1,k+1}^{''}(\Delta_{g(k)}^2-\E(\Delta_0^2))\right)\right\vert
$$
goes to zero as $n$ goes to infinity. Let $M>0$ be fixed. We have $F_2\leq F_2^{'}+F_2^{''}$ where
$$
F_2^{'}=\left\vert\E\left(\frac{1}{n^d}\sum_{k=1}^{n^d}h_{k-1,k+1}^{''}
\left(\Delta_{g(k)}^2-\E_M\left(\Delta_{g(k)}^2\right)\right)\right)\right\vert
$$
and
$$
F_2^{''}=\left\vert\E\left(\frac{1}{n^d}\sum_{k=1}^{n^d}h_{k-1,k+1}^{''}
\left(\E_M\left(\Delta_{g(k)}^2\right)-\E(\Delta_0^2)\right)\right)\right\vert
$$
where we recall the noatation $\E_M\left(\Delta_{g(k)}^2\right)=\E\left(\Delta_{g(k)}^2\vert\F_{V_{g(k)}^M}\right)$. 
The following result is a Serfling type inequality which can be found in \cite{McLeish}.
\begin{lemma}\label{Serfling-inequality}
Let $\U$ and $\V$ be two $\sigma$-algebras and let $X$ be a random variable measurable with respect to $\mathcal{U}$. 
If $1\leq p\leq r\leq\infty$ then 
$$
\|\E(X\vert\V)-\E(X)\|_p\leq2(2^{1/p}+1)\left(\alpha(\U,\V)\right)^{\frac{1}{p}-\frac{1}{r}}\|X\|_r.
$$
\end{lemma}
Applying Lemma $\ref{Serfling-inequality}$ and keeping in mind that $\Delta_0$ is bounded by $4\|\K\|_{\infty}/\sqrt{b_n}$, we derive 
$$
F_2^{''}\leq\|\E_M\left(\Delta_0^2\right)-\E(\Delta_0^2)\|_1\leq\frac{96\|\K\|_\infty^2}{b_n}\,\alpha_{1,\infty}(M)
$$
In the other part,
$$
F_2^{'}=\frac{1}{n^d}\sum_{k=1}^{n^d}\left(J_k^{1}(M)+J_k^2(M)\right)\\
$$
where
$$
J_k^1(M)=\left\vert\E\left(h_{k-1,k+1}^{''M}\circ T^{-g(k)}
\left(\Delta_0^2-\E_M\left(\Delta_0^2\right)\right)\right)\right\vert=0
$$
since $h_{k-1,k+1}^{''M}\circ T^{-g(k)}$ is $\F_{V_0^M}$-measurable and
\begin{align*}
J_k^2(M)&=\left\vert\E\left(\left(h_{k-1,k+1}^{''}\circ T^{-g(k)}-h_{k-1,k+1}^{''M}\circ T^{-g(k)}\right)
\left(\Delta_0^2-\E_M\left(\Delta_0^2\right)\right)\right)\right\vert\\
&\leq\E\left(\left(2\wedge\sum_{\vert i\vert<M}\frac{\vert\Delta_i\vert}{n^{d/2}}\right)\Delta_0^2\right)\\
&\leq\frac{4\|\K\|_{\infty}\,\E(\Delta_0^2)}{(n^db_n)^{1/2}}+\frac{4\|\K\|_{\infty}}{(n^db_n)^{1/2}}\sum_{\substack{\vert i\vert<M \\ i\neq 0}}\E\vert\Delta_i\Delta_0\vert
\quad\textrm{since $\Delta_0\leq\frac{4\|K\|_\infty}{\sqrt{b_n}}$ a.s.}\\
&=O\left(\frac{1}{(n^db_n)^{1/2}}+\frac{M^d\sqrt{b_n}}{n^{d/2}}\right)\qquad\textrm{(by Lemma $\ref{lemme-technique}$)}
\end{align*}
So, putting $M=b_n^{\frac{-1}{2d-1}}$ and keeping in mind that $\sum_{m\geq0}m^{2d-1}\,\alpha_{1,\infty}(m)<+\infty$, we derive
$$
F_2=O\left(M^{2d-1}\,\alpha_{1,\infty}(M)\right)+O\left(\frac{1+b_n^{\frac{d-1}{2d-1}}}{(n^db_n)^{1/2}}\right)=o(1).
$$
The proof of Theorem $\ref{convergence-loi}$ is complete.
\section{Annex}
{\em Proof of Lemma $\ref{lemme-technique}$}. For any $i$ in $\Z^d$ and any $z$ in $\R$, we note $\K_{i}(z)=\K\left(\frac{z-X_i}{b_n}\right)$. So, if 
$s$ and $t$ belongs to $\R$, we have 
$$
\E(Z_0(s)Z_0(t))=\frac{1}{b_n}\bigg(\E\left(\K_{0}(s)\K_0(t)\right)-\E\K_0(s)\E\K_0(t)\bigg)
$$
and
$$
\lim_{n\to+\infty}\frac{1}{b_n}\E\left(\K_0(s)\K_0(t)\right)=\lim_{n\to+\infty}\int_{\R}\K \left(v\right)\K \left(v+\frac{t-s}{b_n}\right)f(s-vb_n)dv
=\delta_{st}\,f(s)\,\sigma^2
$$
where $\delta_{st}=1$ if $s=t$ and $\delta_{st}=0$ if $s\neq t$. We have also
$$
\lim_{n\to+\infty}\frac{1}{b_n}\E\K_0(s)\E\K_0(t)=\lim_{n\to+\infty}b_n\int_{\R}\K (v)f(s-vb_n)dv\int_{\R}\K (w)f(t-wb_n)dw=0.
$$
So, we obtain
$$
\E (\Delta_0^2)=\lambda_1^2\E (Z_0^2(x))+\lambda_2^2\E (Z_0^2(y))+2\lambda_1\lambda_2\E(Z_0(x)Z_0(y))\converge{n}{+\infty}{ }\eta.
$$
Let $i\neq 0$ be fixed in $\Z^d$. We have 
\begin{equation}\label{delta_O_delta_i}
\E\vert\Delta_0\Delta_i\vert\leq\lambda_1^2\E\vert Z_0(x)Z_i(x)\vert+\lambda_2^2\E\vert Z_0(y)Z_i(y)\vert
+\lambda_1\lambda_2\E\vert Z_0(x)Z_i(y)\vert+\lambda_1\lambda_2\E\vert Z_0(y)Z_i(x)\vert.
\end{equation}
For any $s$ and $t$ in $\R$, 
$$
\E\vert Z_0(s)Z_i(t)\vert\leq\frac{1}{b_n}\E\big\vert \K_0(s)\K_i(t)\big\vert+\frac{1}{b_n}\E\big\vert\K_0(s)\big\vert\,\E\big\vert \K_0(t)\big\vert.
$$
Moreover, using Assumptions $\textbf{(A2)}$ and $\textbf{(A3)}$, we have
$$
\frac{1}{b_n}\E\big\vert\K_0(s)\big\vert\,\E\big\vert \K_0(t)\big\vert=b_n\int_{\R}\vert \K (u)\vert f(s-ub_n)du\int_{\R}\vert \K (v)\vert f(t-vb_n)dv=O(b_n)
$$
and
$$
\frac{1}{b_n}\E\big\vert \K_0(s)\K_i(t)\big\vert=b_n\iint_{\R^2}\big\vert\K\left(w_1\right)\K \left(w_2\right)\big\vert f_{0,i}(s-w_1b_n,t-w_2b_n) dw_1dw_2=O(b_n).
$$
So, we obtain for any $s$ and $t$ in $\R$
\begin{equation}\label{Z_0_Z_i}
\E\vert Z_0(s)Z_i(t)\vert=O(b_n).
\end{equation}
The proof of Lemma $\ref{lemme-technique}$ is completed by combining ($\ref{delta_O_delta_i}$) and ($\ref{Z_0_Z_i}$).\\
\\
{\em Proof of Lemma $\ref{lemme-conditions-on-m_n}$}. First, $m_n^d$ goes to infinity since $b_n$ goes to zero and $m_n\geq\left[b_n^{-\frac{1}{2d}}\right]$. 
For any positive integer $m$, we consider
$$
\psi(m)=\sum_{\vert i\vert>m}\vert i\vert^d\,\alpha_{1,\infty}(\vert i\vert).
$$
Since the mixing condition ($\ref{mixing-condition}$) is equivalent 
to $\sum_{k\in\Z^d}\vert k\vert^d\,\alpha_{1,\infty}(\vert k\vert)<\infty$, we know that $\psi(m)$ converges to zero as $m$ goes to infinity. Moreover, we have
$$
m_n^db_n\leq\max\left\{\sqrt{b_n},\sqrt{\psi\left(\left[b_n^{-\frac{1}{2d}}\right]\right)}+2^db_n\right\}\converge{n}{+\infty}{ }0.
$$
We have also
$$
m_n^d\geq\frac{1}{b_n}\sqrt{\psi\left(\left[b_n^{-\frac{1}{2d}}\right]\right)}\geq\frac{1}{b_n}\sqrt{\psi\left(m_n\right)}
\quad\textrm{since $\left[b_n^{-\frac{1}{2d}}\right]\leq m_n$}.
$$
Finally, we obtain
$$
\frac{1}{m_n^db_n}\sum_{\vert i\vert>m_n}\vert i\vert^d\,\alpha_{1,\infty}(\vert i\vert)\leq\sqrt{\psi(m_n)}\converge{n}{+\infty}{ }0.
$$
The proof of Lemma $\ref{lemme-conditions-on-m_n}$ is complete.
\bibliographystyle{plain}
\bibliography{xbib}
\end{document}